\documentclass{ifacconf}

\usepackage{graphicx}      % include this line if your document contains figures
\usepackage{mathptmx}
\usepackage{url}
\usepackage[english]{babel}
\usepackage{natbib}        % required for bibliography
\usepackage{array}
\usepackage{pagecolor}
\usepackage{amsmath,amsfonts,latexsym,amscd,amssymb,theorem}
\usepackage{epsfig}
\usepackage{color}
\usepackage{amsmath}
\usepackage{longtable}
\usepackage{booktabs}
\usepackage{bbm}
\usepackage[T1]{fontenc}
\usepackage[utf8]{inputenc}
\usepackage{esvect}
 \usepackage{float}
\usepackage{multirow}
\usepackage{kpfonts}
\usepackage{tikz}
\usepackage{babelbib}
\usepackage{algorithm}
\usepackage{algpseudocode}

\DeclareMathOperator*{\argmin}{arg\,min}

\def\ls{\textcolor{black}}

\begin{document}
\begin{frontmatter}

\title{Optimizing semilinear representations for State-dependent Riccati Equation-based feedback control} 

\thanks[footnoteinfo]{This work was supported
	by the UK Engineering and Physical Sciences Research Council (EPSRC) grant
	EP/V04771X/1. DK was additionally supported by EPSRC grants EP/T024429/1, and EP/V025899/1.}

\author[First]{S. Dolgov} 
\author[Second]{D. Kalise} 
\author[Second]{L. Saluzzi}

\address[First]{Department of Mathematical Sciences, University of Bath, North Rd, BA2 7AY Bath, United Kingdom, (e-mail: {sd901}@bath.ac.uk).}
\address[Second]{Department of Mathematics, Imperial College London, South Kensington Campus, SW7 2AZ London, United Kingdom, (e-mail: dkaliseb@imperial.ac.uk, lsaluzzi@ic.ac.uk).}
%\address[Third]{Department of Mathematical Sciences, University of Bath, North Rd, BA2 7AY Bath, United Kingdom, (e-mail: @bath.ac.uk)}

\begin{abstract}     
           % Abstract of not more than 250 words.
$\;$ An optimized variant of the State Dependent Riccati Equations (SDREs) approach for nonlinear optimal feedback stabilization is presented. The proposed method is based on the construction of equivalent semilinear representations associated to the dynamics and their affine combination. The optimal combination is chosen to minimize the discrepancy between the SDRE control and the optimal feedback law stemming from the solution of the corresponding Hamilton Jacobi Bellman (HJB) equation. Numerical experiments assess effectiveness of the method in terms of stability of the closed-loop with near-to-optimal performance.
\end{abstract}


\end{frontmatter}
%===============================================================================

\section{Introduction}
Feedback control laws are fundamental for asymptotic stabilization of nonlinear dynamics. Their synthesis can be performed in an ad-hoc manner leading to suboptimal controllers, or using optimal control methods, in which case the final outcome is an optimal stabilizing feedback law.   The natural framework for the computation of optimal stabilizing controllers is the use of dynamic programming techniques, leading ultimately to the solution of a Hamilton-Jacobi-Bellman partial differential equation (HJB PDE). This poses a formidable computational challenge, as the HJB PDE is a nonlinear first order PDE cast over the state space of the system, with a dimension which can be arbitrarily high. Over the last years, the solution of large-scale optimal feedback synthesis problems has witnessed a tremendous development, in parallel with the development of sophisticated techniques for high-dimensional problems, such as representation formulas for the HJB PDE \cite{chow1}, tensor decomposition techniques \cite{tensor4}, tree-structured algorithms \cite{AS20} and data-driven methods \cite{ml1,AKK21,nakazim}, among others. An alternative approach, which can be interpreted as a relaxed dynamic programming, is provided by Nonlinear Model Predictive Control (NMPC) \cite{GR2008}. Here, a control signal is optimized over a prediction horizon, and updated as the dynamics evolve. This mechanism generates a feedback law in the sense that a sequential re-computation of the control signal accounts for perturbations along the trajectory. A proper selection of the prediction horizon, as well as the running costs to be optimized, can guarantee asymptotic stabilization of the closed-loop \ls{(see for instance Theorem 4.8 in \cite{larsnmpc})}. 

In this work, we approach the optimal feedback stabilization problem by resorting to an alternative which combines dynamic programming and NMPC element, known as the State-Dependent Riccati Equation (SDRE) approach \cite{sdrefirst1,sdrefirst2}. The SDRE feedback law is constructed upon an approximation of the HJB PDE, thus constituting a suboptimal feedback control. However, under suitable stabilizability assumptions, it can generate a locally asymptotically stabilizing feedback law. Formally, the SDRE synthesis is based around the sequential solution of algebraic Riccati equations along a trajectory, where the operators involved are state-dependent functions which are frozen at the current state of the system. It is thus, an implementation which is reminiscent of the NMPC closed-loop. The SDRE methodology has been extensively studied in \cite{WANG98,CIMEN,BTB00,BLT07} and continues to be an active subject of research due to its simplicity and effectiveness \cite{HK18,heiland,Astolfi2020,ABK21}. \ls{A similar approach is represented by the extension of Al’brekht’s Method considering higher order Taylor expansions, as discussed in \cite{krener2020brekht}.}

Here, we study an issue which has been extensively discussed in the SDRE literature. The SDRE feedback law requires the nonlinear dynamics to be expressed in semilinear form, that is $f(x)=A(x)x$, for $A(x)$ a state-dependent matrix. However, this semilinear representation is non-unique, and its selection can have an impact on the SDRE feedback law. In the worst case scenario, a poor selection of a representation can hinder the stability of the resulting-closed loop. However, as recently shown in \cite{Astolfi2020}, it can be favourably used to pick a representation which minimizes the misfit between the SDRE and HJB feedback laws, this latter being the true object we wish to approximate. Following this vein, here we present a systematic procedure to generate a family of semilinear representations, for which we consider an affine combination where the coefficients are optimized to minimize the discrepancy with respect to the HJB feedback law. We present nonlinear numerical tests where we show that this approach offers a suitable alternative to improve the performance of the SDRE controller, ensuring stabilization and close-to-optimal performance.

\section{State Dependent Riccati Equation}
Let us consider a control-affine dynamical system driven by the following system of ordinary differential equations in $\mathbb{R}^d$:
\begin{equation}
\left\{ \begin{array}{l}
\dot{y}(t)=f(y(t)) + B(y(t))u(t), \;\; t\in(0,+\infty),\\
y(0)=y^0\in\mathbb{R}^d,
\end{array} \right.
\label{dyn}
\end{equation}
where the control signal $u \in \mathcal{U}=\{u: [0, +\infty) \rightarrow \mathbb{R}^m$, measurable $\}$. We will assume that $f(\underline{0})=\underline{0}$, $i.e.$ the origin is an equilibrium for the dynamical system \eqref{dyn}. 
We are interested in asymptotic stabilization to the origin by minimizing the following cost functional
$$
J(u;x):=\frac{1}{2}\int_0^{+\infty} y(t)^{\top} Q y(t)+ u^{\top}(t) u(t)\,dt\,,
$$
subject to \eqref{dyn}, where $Q \in \mathbb{R}^{d \times d}$ is a positive semidefinite matrix. This problem is solved using dynamic programming. For this, we define the value function of the problem
$$
V(x)= \inf_{u \in \mathcal{U}} J(u(\cdot,x)).
$$
It is well known that the value function satisfies the following nonlinear Hamilton-Jacobi-Bellman PDE:
\begin{equation}\label{HJB}
\min\limits_{u\in \mathbb{R}^m }\left\{(f(x)+B(x)u)^{\top} \nabla V(x) +\frac{1}{2} x^{\top} Q x+ \frac{1}{2}u^{\top}  u \right\}=0.
\end{equation}
After solving equation \eqref{HJB}, the optimal feedback control is expressed in feedback form as
\begin{equation}
u^*(x)=-B(x)^{\top} \nabla V(x)\,.
\label{opt_con}
\end{equation}
Equation \eqref{HJB} is a nonlinear first order PDE cast over the state space of the dynamics \eqref{dyn}, hence sharing its dimension and eventually suffering from the \textsl{curse of dimensionality}. For this reason, we avoid the numerical approximation of the HJB PDE by traditional discretization methods and instead we resort to the use of the State-dependent Riccati Equation (SDRE) approach as a way to generate a suboptimal approximation of the optimal stabilizing law. For this, we write the dynamics \eqref{dyn} in semilinear form $f(x)=A(x)x$ leading to
\begin{equation}
\dot{y}(t)=A(y(t))y(t) + B(y(t))u(t), \;\; t\in(0,+\infty)\,.
\end{equation}
Assuming that for each point $x \in \mathbb{R}^d$ the couple $(A(x),B(x))$ is stabilizable \ls{and detectable}, it is possible to construct a locally asymptotically stabilizing feedback control by solving the State-dependent Riccati equation:
\begin{equation}
A^{\top}(x) \Pi(x) + \Pi(x) A(x) - \Pi(x) W(x) \Pi(x)+Q=0,
\label{sdre}
\end{equation}
where $W(x)=B(x) B^{\top}(x)$.
In the vast majority of cases of interest, the SDRE \eqref{sdre} cannot be solved analytically, and the feedback law is implemented in a model predictive control fashion, by measuring the current state $x$ of the system, solving \eqref{sdre} for a fixed state, applying the resulting feedback law, and evolving the trajectory.
From the SDRE, an approximate value function $\widetilde{V}(x)$ can be obtained as
\begin{equation}
\widetilde{V}(x)=\frac{1}{2}x^{\top} \Pi(x) x.
\label{vf2}
\end{equation}

As discussed in \cite{Astolfi2020}, the residual in approximating the value function $V(x)$, which solves \eqref{HJB}, using the expression \eqref{vf2} is defined as
\begin{equation}
E(x)=\varphi(x)\left(2\left[A(x)-W(x) \Pi(x)\right]x -W(x) \varphi(x)^T\right),
\label{error}
\end{equation}
where
\begin{equation}
(\varphi(x))_k= \frac{1}{2}\sum_i \sum_j x_i x_j \partial_{x_k} \Pi_{i,j}(x).
\label{phi}
\end{equation}
The term $\partial_{x_k} \Pi_{i,j}(x)$ denotes the partial derivative with respect to $x_k$ of the component  $(i,j)$ of the matrix $\Pi(x)$ and can be computed deriving equation \eqref{sdre} with respect to $x_k$, obtaining the following Lyapunov equation:
\begin{equation}
\partial_{x_k}  \Pi(x) \left(A(x)-W \Pi(x) \right) + \left(A(x)^{\top}-\Pi(x) W \right) \partial_{x_k}  \Pi(x) + \mathcal{Q}_k=0,
\label{Px}
\end{equation}
where
$$
\mathcal{Q}_k=\partial_{x_k} A(x)^{\top} \Pi(x) + \Pi(x) \partial_{x_k} A(x)-\Pi(x)\partial_{x_k}W(x)  \Pi(x) .
$$

Following \eqref{opt_con} and \eqref{vf2}, the feedback map constructed by this procedure has the following form:
\begin{equation}
u^*(x)=-B(x)^{\top} \left( \Pi(x)x + \varphi(x) \right).
\label{opt_con2}
\end{equation}
\section{Design of semilinearization for SDRE}
The SDRE is constructed upon a semilinear representation of the dynamics $f(x)=A(x)x$, however, this representation is non-unique and this can have an impact in the control law. In the following, we present a synthesis method where this non-uniqueness property is used to our benefit, by choosing a representation where the misfit between the SDRE and the HJB value function is minimized.
Let us suppose we have $N+1$ possible semilinear representations of the nonlinear vector field $f(x)$, $i.e.$ 
$$
f(x)=A_i(x)x,\quad i\in \{0,\ldots, N\}\,.
$$
It is easy to see that the equality above still holds for an affine combination, $i.e.$
$$
f(x)= \sum_{i=0}^N \alpha_i A_i(x) x,
$$
with
$$
\alpha=(\alpha_i)_i \in \mathcal{A}_N=\left\{ \alpha \in \mathbb{R}^{N+1}:  \sum_{i=0}^N \alpha_i=1 \right\} \,.
$$
Let us define the affine combination of the matrices $\{A_i(x)\}_i$ as 
$$
\mathcal{A}(x,\alpha)=\sum_{i=0}^N \alpha_i A_i(x).
$$
\ls{We suppose that exists at least one combination $\tilde{\alpha}=(\tilde{\alpha}_i)_i$ such that the couple $(\mathcal{A}(x,\tilde{\alpha}),B(x))$ is  stabilizable and detectable.}
Exploiting the equality constraint, the problem can be formulated in a equivalent way as
$$
\mathcal{A}(x,\alpha)=\sum_{i=0}^{N-1} \alpha_i A_i(x)+\left(1-\sum_{i=0}^{N-1} \alpha_i \right) A_N(x), \quad \alpha \in \mathbb{R}^N.
$$

Fixing $x \in \mathbb{R}^d$ and $\alpha \in  \mathbb{R}^N $, the solution of the SDRE $\Pi(x,\alpha)$, its derivatives $\{\partial_{x_k} \Pi(x,\alpha)\}_k$ and the residual $E(x,\alpha)$ can be computed following the steps described in the previous section. Having constructed a family $\{A_i(x)\}_i$, our aim is to solve the following minimization problem for every point $x$ along a given trajectory:
\begin{equation}
\label{min_E}
E^*(x)=\min_{\alpha \in \mathbb{R}^N} E(x,\alpha)^2 .
\end{equation}
The convergence of $E^*(x)$ to zero will imply the convergence of the approximated value function $\widetilde{V}(x)$ to $V(x)$, retrieving an optimal feedback law.
In Algorithm \ref{alg_1} we sketch the method to compute the optimal trajectory based on the minimization \eqref{min_E}.

\begin{algorithm}[H]
%\captionsetup[algorithm]{name=TSA}
\caption{Optimal trajectory based on optimized semilinear representations}
\begin{algorithmic}[1]
\State Choose a collection $\{A_i(x)\}_{i=0}^N$, an initial point $y^0$, the number of time steps $n_t$, a tolerance $tol$ and a initial combination $\alpha^*_{-1}$
\State Define $\mathcal{A}(x,\alpha)=\sum_{i=0}^N \alpha_i A_i(x) $
\For{$i=0,\ldots, n_t-1$}
\If{$E(y^{i},\alpha^*_{i-1})^2 \le tol$\label{if1}}
\State{$\alpha^*_i :=\alpha^*_{i-1}$ \label{if2}}
\Else
\State{$\alpha_i^* \in \argmin_{\alpha} E(y^{i},\alpha)^2$ \label{min}}
\EndIf
\State{Define $A(y^{i})=\mathcal{A}(y^{i},\alpha_i^*)$}
\State{Solve \eqref{sdre} obtaining $\Pi(y^{i})$}
\State{Solve \eqref{Px} obtaining $\{\partial_{k}\Pi(y^{i})\}_k$ and $\varphi(y^{i})$}
\State{Compute $u^*(y^i)$ via \eqref{opt_con2}}
\State{Evolve to $y^{i+1}$}
\EndFor
\end{algorithmic}
\label{alg_1}
\end{algorithm}

Line \ref{if1}-\ref{if2} of Algorithm \ref{alg_1} are introduced to avoid unnecessary computations, keeping the same minimizer if the residual computed on the new point of the trajectory stays below a threshold $tol$.
Line \ref{min} requires the resolution of a minimization problem in dimension $N$ and it does not require any bound, since the equality $\mathcal{A}(x,\alpha)x=f(x)$ holds for every $\alpha \in \mathbb{R}^N$.  However, $E(x,\alpha)^2$ is non-convex and may have different local minima. In this case one may consider different random initializations or a global optimization solver to achieve a better accuracy. Having fixed the optimal parameter $\alpha^*$, for each time step we will need to solve a Riccati equation \eqref{sdre} and $d$ Lyapunov equations \eqref{Px} to assemble the optimal control.

A relevant issue in this formulation is the systematic generation of a family of semilinear representations of the nonlinear dynamics. First, we fix a first matrix $A_0(x)$ such that $f(x)=A_0(x)x$. Starting from $A_0(x)$, it is possible to construct an equivalent semilinearization. Here we propose the following procedure: \ls{choosing a row index $i_1 \in \{1,\ldots,d\}$, two column indices $j_1 \in \{1,\ldots, d-1\}$ and $j_2 \in \{j_1+1,\ldots, d\}$ and an arbitrary scalar function $h(x) \in \mathcal{C}$}, it is possible to check that the following matrix
$$
[\tilde{A}(x)]_{i,j}= \begin{cases} [A_0(x)]_{i,j}+h(x) x_{j_2} & (i,j)=(i_1,j_1),  \\
[A_0(x)]_{i,j}-h(x) x_{j_1}, & (i,j)=(i_1,j_2), \\
[A_0(x)]_{i,j} & otherwise
\end{cases}
$$
still represents a semilinearization for the vector field $f$.
\ls{In this way it is possible to construct up to  $N=\frac{d^2(d-1) |\mathcal{C}|}{2}$ matrices using this procedure, where $|\mathcal{C}|$ represents the cardinality of the set $\mathcal{C}$.}

Solving the entire $M=N$-dimensional optimization problem \eqref{min_E} may be difficult for large \ls{$M=\mathcal{O}(d^3|\mathcal{C}|)$.}
Instead, we can restrict $\alpha$ to the form $(0,\ldots,0,\alpha_i,0,\ldots,0)$ and solve up to $N$ one-dimensional optimization problems over $\alpha_i$ until the residual $E(x,\alpha)^2$ is below the desired threshold. \ls{The choice of the number of semilinear forms is arbitrary, according to the desired efficiency and computational cost.}

In the numerical tests we will consider $h(x)=c \in \mathbb{R}$ and in particular $c \in \{-1,1 \}$, since we are interested in linear perturbation of the matrix $A_0(x)$.
%To accelerate the computations,
%instead of solving the full $N$-dimensional optimization problem \eqref{min_E} we can consider a simplified method.
% As we will see in the next section, in the presented cases it will be sufficient to consider a single permutation at a each time, obtaining a one-dimensional minimization problem. We will vary among the different single permutations until the error is below a fixed threshold.
\section{Numerical tests}
In this section we assess the performance of the proposed methodology in two nonlinear tests. The minimization \eqref{min_E} is solved using the Matlab function \texttt{fminunc} \ls{with initial guess $\alpha^0_i=1$ $\forall i=0,\ldots,N$} and the solution for the differential equation \eqref{dyn} is obtained using the Matlab function \texttt{ode45}. The numerical simulations reported in this paper are performed on a Dell XPS 13 with Intel Core i7, 2.8GHz and 16GB RAM. The codes are written in Matlab R2022a.
\subsection{Lorentz system}

In the first example we consider the Lorentz system
\begin{equation*}
\begin{cases}
\dot{x}= \sigma (y-x), \\
\dot{y}= x(\rho-z)-y+u, \\
\dot{z}=xy-\beta z
\end{cases}
\end{equation*}
and the cost functional
\begin{equation*}
J= \frac12\int_0^\infty 100( |x(s)|^2+|y(s)|^2+|z(s)|^2 )+ |u(s)|^2 \, ds.
\end{equation*}
We fix $\sigma=10$, $\beta=8/3$, $\rho=2$ and $(x_0,y_0,z_0)=(-1,-1,-1)$.
We fix the following matrix for the initial semilinear representation:
$$
A_0(x,y,z)=\begin{bmatrix} -\sigma & \sigma & 0 \\ \rho-z  &  -1 & 0 \\ y& 0 &-\beta \end{bmatrix}.
$$

\begin{figure}[htbp]
         \centering
\includegraphics[scale=0.4]{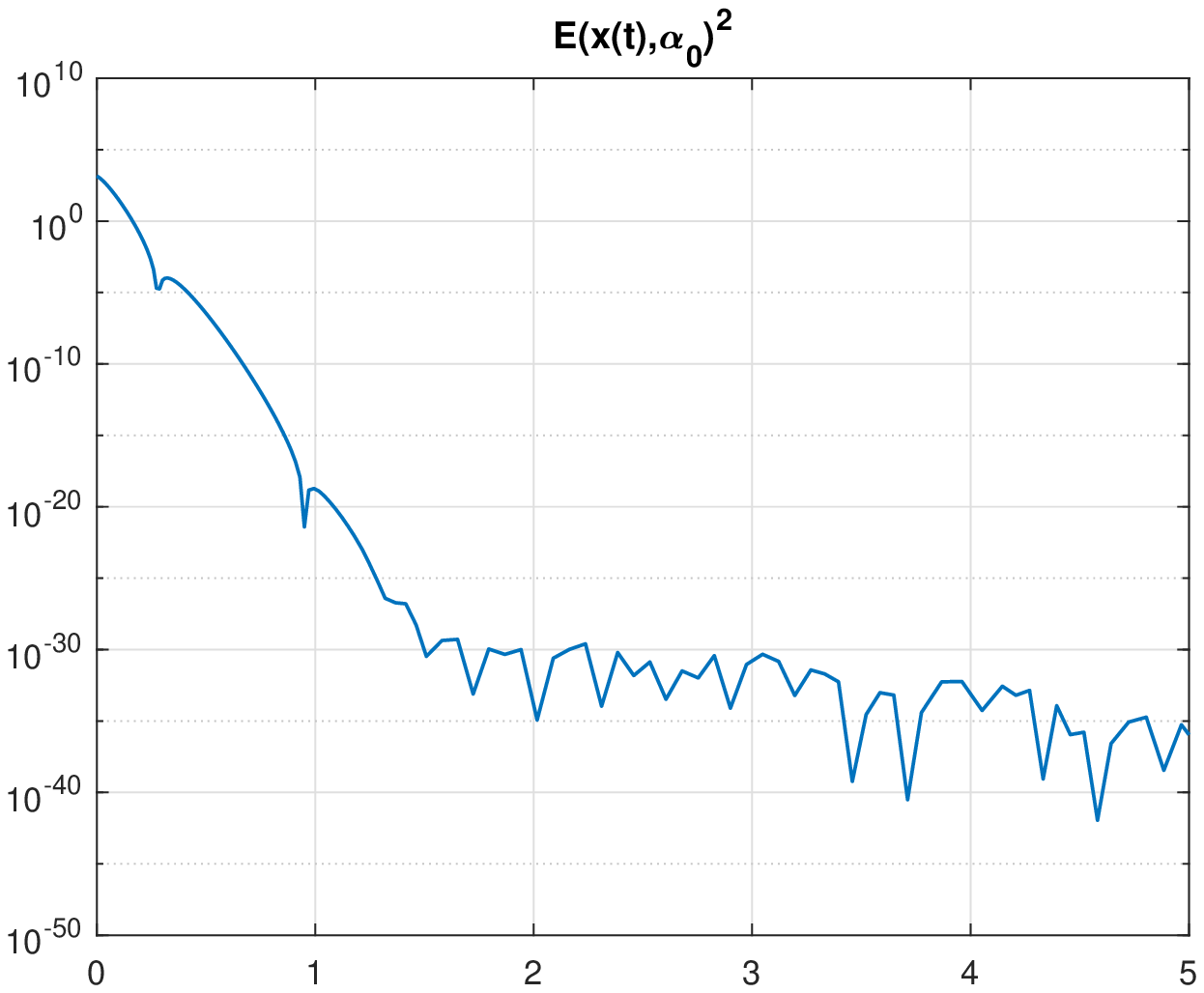}
\includegraphics[scale=0.4]{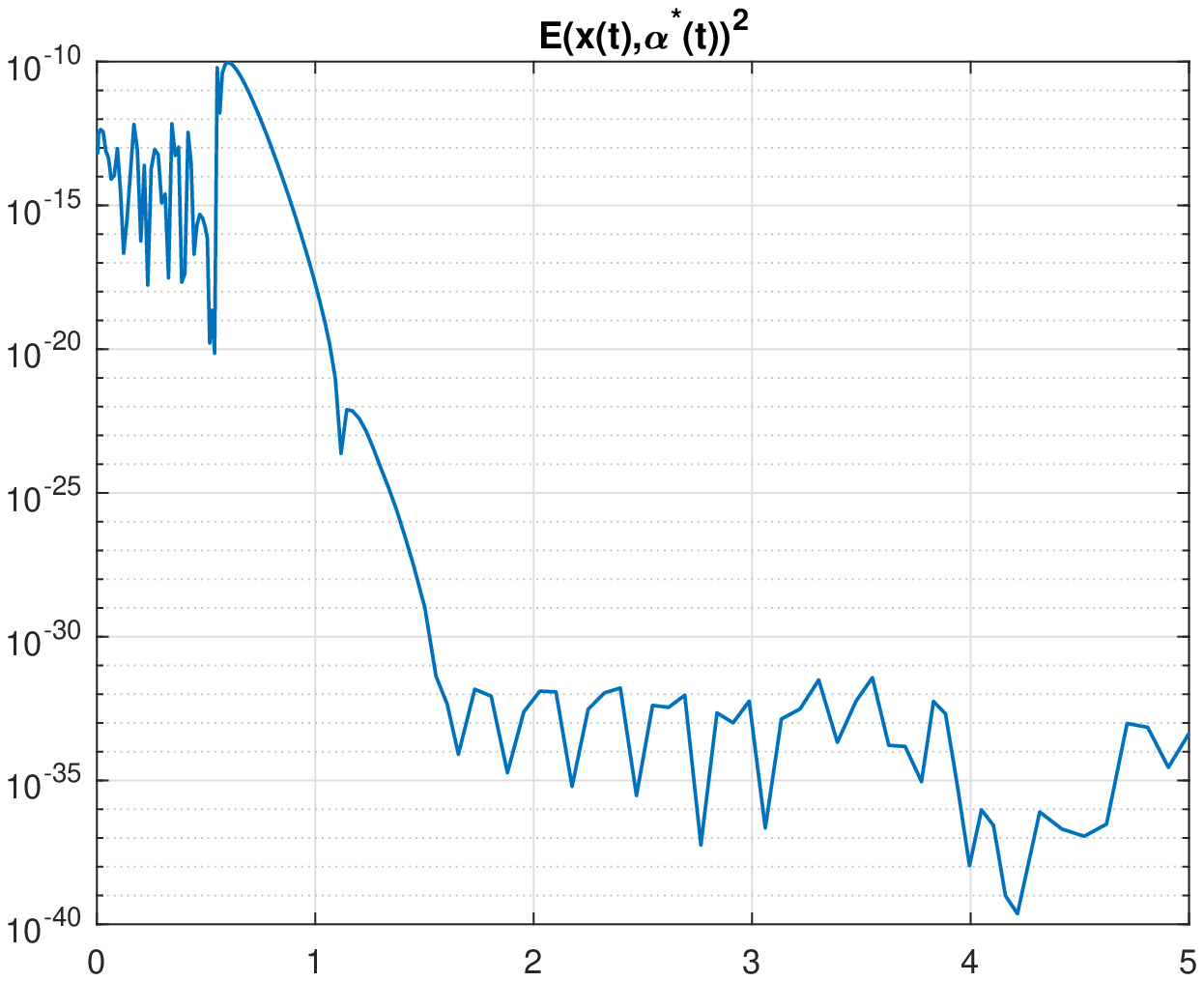}

                \caption{Test 1: Residual along the optimal trajectory with $\alpha_0$ (top) and with the optimal $\alpha^*(t)$ (bottom) in logarithmic scale.}
       \label{fig:case1}
	\end{figure}

\begin{figure}[htbp]
         \centering
       \includegraphics[scale=0.4]{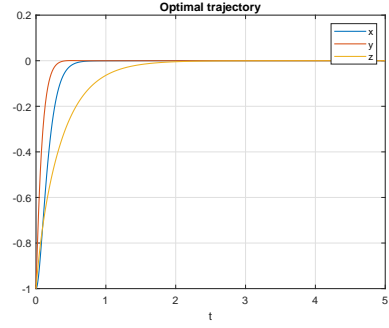} 
                \caption{Test 1:  Optimal trajectory with $\alpha^*(t)$ and initial condition $(-1,-1,-1)$}
       \label{fig:case2}
	\end{figure}

We notice that the matrix $A_0$ depends linearly on the variables, hence some of the possible linear perturbations discussed previously will cancel some terms. We will compare Algorithm \ref{alg_1} applied with the fixed semilinear representation $A_0(x)$ ($i.e.$ choosing $\alpha=(1,0,\ldots,0)$) against the minimization over affine combinations. We fix $tol=10^{-12}$. In Table \ref{Table_err} we compare these two cases in terms of total cost $J$ and total residual $\int E(y(s),\alpha(s))^2 ds$ computed along the optimal trajectory. The introduction of the minimization leads to a total residual of order $10^{-12}$, while the fixed semilinear representation is characterized by a much higher residual. This is also reflected in the computation of the total cost, obtaining an improvement of almost $10\%$ with respect to the fixed case. In Figure \ref{fig:case1} we show the behaviour of the function $E(x(t),\alpha(t))$ along the dynamics for the fixed case (top panel) and for the optimal one (bottom panel). We notice that by choosing the optimal $\alpha^*$, the order of the residual is always lower than $10^{-10}$, while we can see in both cases that as the dynamics gets closer to the origin, the residual reaches a very low value. Figure \ref{fig:case2} displays the optimal trajectory using the optimal affine combination $\alpha^*$, showing the convergence of the entire system to the equilibrium. %In the central panel we display the behaviour of the optimal parameter $\alpha^*(t)$ and in the bottom panel the index of the selected permutation, characterized by a residual under the threshold $tol$. In this case we are considering all the possible permutations with $h(x)=c \in \{-1, 1\}$, obtaining $M=18$ semilinear represenations. We notice that starting from time $t \approx 0.8$, the minimization is always solved by the first permutation and the optimal solution is  given by $\alpha^*=1$.
%In the central panel and in the bottom panel of Figure \ref{fig:case2} we show the profile of the error $E(x,\alpha)$ for two different settings. It is possible to notice the presence of multiple local minima distributed on a 1 dimensional curve. This minimization problem requires different initializations and their comparisons in order to obtain an accurate resolution.

\begin{table}[htbp]			
			
		\centering
		\begin{tabular}{c|cc}
					\toprule
					Semilinear form &Total cost & Total res  \\
					\midrule
$A_0(x)$ &  5.79 &  45.8\\
$\mathcal{A}(x,\alpha^*)$ & \textbf{5.27} &   \textbf{ 7.6e-12}\\

			\bottomrule		
		\end{tabular}	
	\caption{Test 1: Comparison between the fixed and the optimal case in terms of total cost and total residual. \label{Table_err}}
	\end{table}

%Finally we note that the computation of the integrals 
%$$
%\int_0^T\alpha_i(s) \,ds \quad \forall i \in \{0,\ldots,N\}
%$$
%gives an idea of the contributions of the matrices $\{A_i(x)\}_i$ in the minimization \eqref{min_E}. In this example it turns to be
%$$
%A(x,y,z)=\begin{bmatrix} -\sigma+y & \sigma-x & 0 \\ \rho-z  &  -1 & 0 \\ y& 0 &-\beta \end{bmatrix} .
%$$ 
%		  \begin{figure}[htbp]
%         \centering
%       \includegraphics[scale=0.4]{pics/opt_traj.eps} 
%\includegraphics[scale=0.4]{pics/mesh_err.eps}
%\includegraphics[scale=0.4]{pics/mesh_err_3d2.eps}
%                \caption{Test 1:  Optimal trajectory with $\alpha^*(t)$ (top), the error in logarithmic scale with $N=2$ and $x=x_0$ (central) and with $N=3$, $x=x_0$ and $\alpha_2$ fixed (bottom).}
%       \label{fig:case2}
%	\end{figure}

\subsection{The inverted pendulum}
The second example deals with the optimal control of the inverted pendulum. \ls{We denote by $x_1$ and $x_3$ as the position of the cart and its velocity, while by $x_2$ and $x_4$ the tilt angle and the corresponding velocity.}
Defining $x=(x_1,x_2,x_3,x_4) \in \mathbb{R}^4$, a possible semilinearization of the dynamics is the following
$$
\dot{x}=A_0(x)x+g(x)x,
$$
where
$$
A_0(x)=\begin{bmatrix} 0 & 0 & 1 &0 \\ 0 & 0 & 0 &1 \\ 0 & a_{3,2}(x) & 0 & 0\\ 0 & a_{4,2}(x) & 0 &a_{4,4}(x) \end{bmatrix} ,
$$
$$
a_{3,2}(x)=m\, s(x)(\ell x_4-gcos(x_2))/c(x),
$$
$$
a_{4,2}=s(x)(M+m)g/(\ell c(x)),
$$
$$
a_{4,4}= -m\,x_4sin(x_2)cos(x_2)/c(x),
$$
$$
g(x)=\begin{bmatrix} 0 \\ 0 \\ 1/c(x) \\ -cos(x_2)/(\ell c(x)) \end{bmatrix},
$$
$$
s(x)=sin(x_2)/x_2, \; c(x)=M+m\,sin(x_2)^2.
$$
We fix $M=0.5$, $m=0.45$, $\ell=0.5$ and $g=9.81$.

\ls{We are interested in the minimization of the following cost functional
\begin{equation*}
J= \frac12\int_0^\infty x^\top(s) Q x(s) + |u(s)|^2 \, ds
\end{equation*}
where $Q=diag(1,10,0.1,0.1)$.
}
%The requirement of the semilinear form leads to the presence of the term $sin(x_2)/x_2$, which is bounded as $x_2$ tends to zero.

\begin{figure}[htbp]
         \centering
       \includegraphics[scale=0.4]{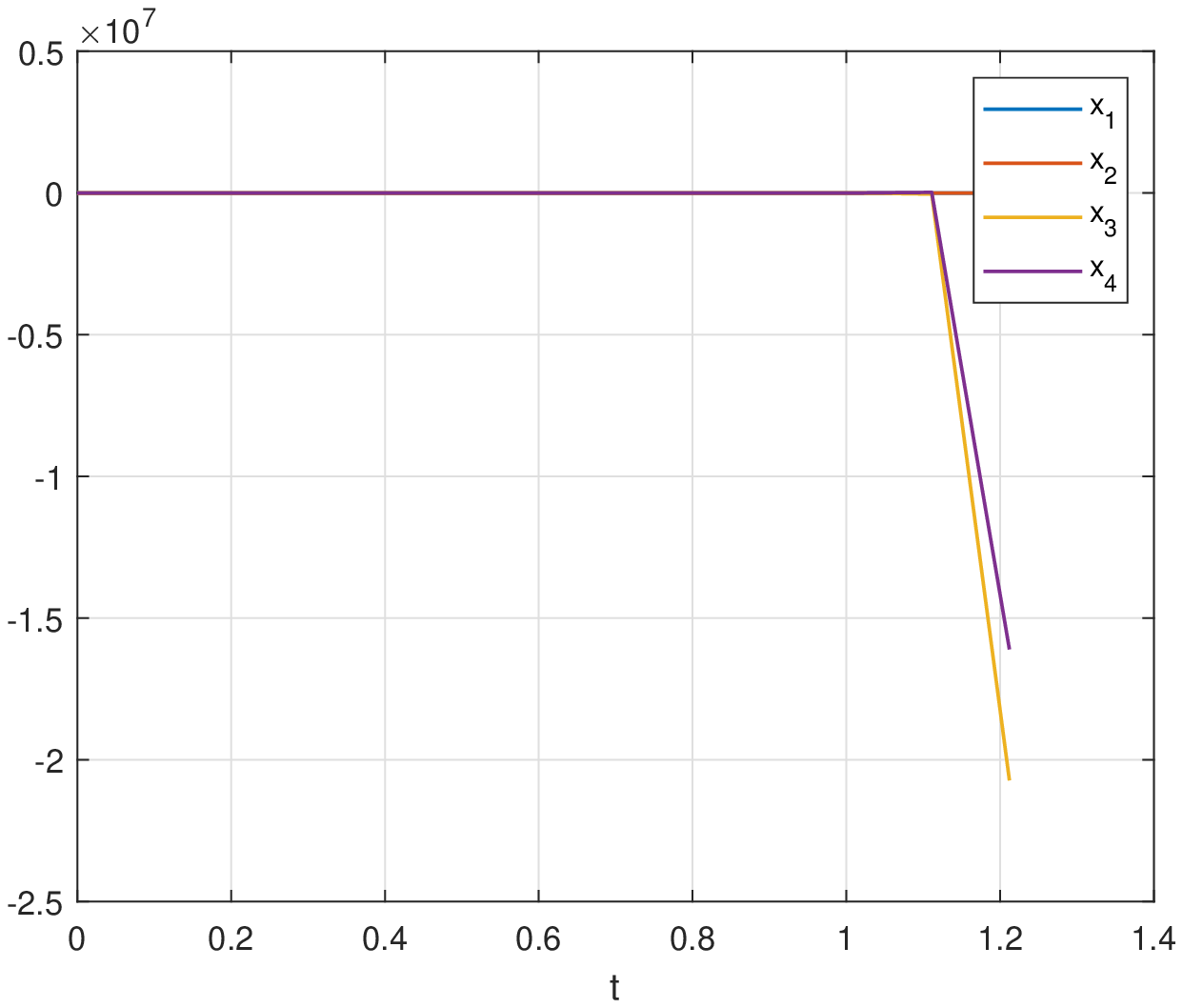} 
\includegraphics[scale=0.4]{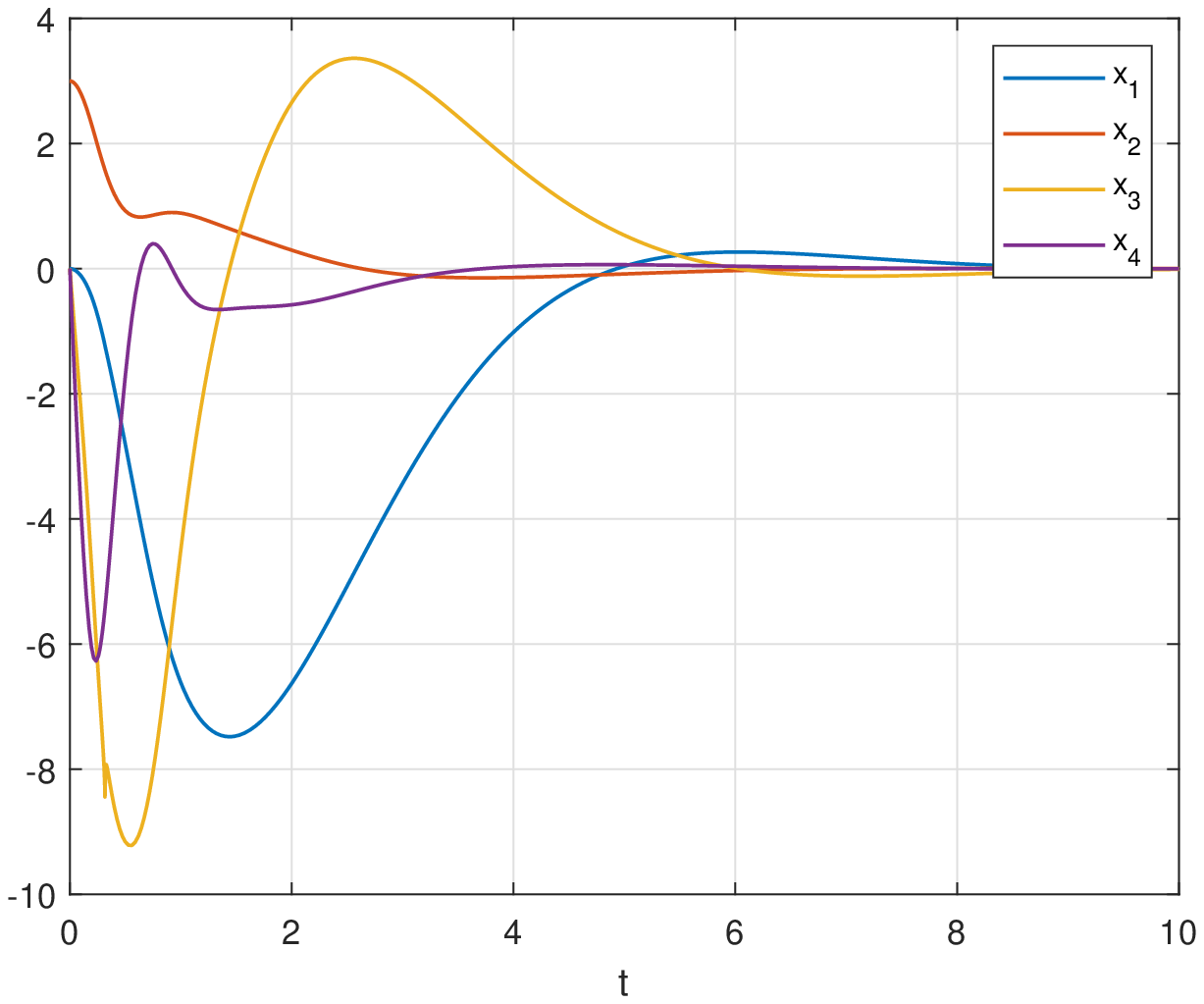}
                \caption{Test 2:  Optimal trajectories with initial condition $(0,3,0,0)$ and different semilinear representations: the case with $A_0(x)$ (top) and the one with $\tilde{A}(x)$ (bottom).}
       \label{fig:case3}
	\end{figure}

We are going to compare again the performance of the algorithm in absence and in presence of the optimal affine combination. We will fix $tol=10^{-9}$. First of all, let us fix the initial condition $(0,3,0,0)$. In this case the construction of the optimal trajectory using $A_0(x)$ fails, since the dynamics passes through points in which the couple $(A_0(x), g(x))$ is not stabilizable and it is not possible to solve the Riccati equation \eqref{sdre}. However, using instead the semilinear representation
$$
[\tilde{A}(x)]_{i,j}= \begin{cases} [A_0(x)]_{i,j}- x_{4} & (i,j)=(2,2),  \\
[A_0(x)]_{i,j}+ x_{2}, & (i,j)=(2,4), \\
[A_0(x)]_{i,j} & otherwise
\end{cases}
$$

the couple $(\tilde{A}(x), g(x))$ becomes stabilizable for all the points along the trajectory. The plots of the optimal trajectories in these two cases are shown in Figure \ref{fig:case3}. In the first case the dynamical trajectory diverges at time $t=1.2$, while the latter case is able to reach the equilibrium. This illustrates that the proposed technique is able to find a stabilizing linear perturbation.

\begin{figure}[htbp]
         \centering
       \includegraphics[scale=0.4]{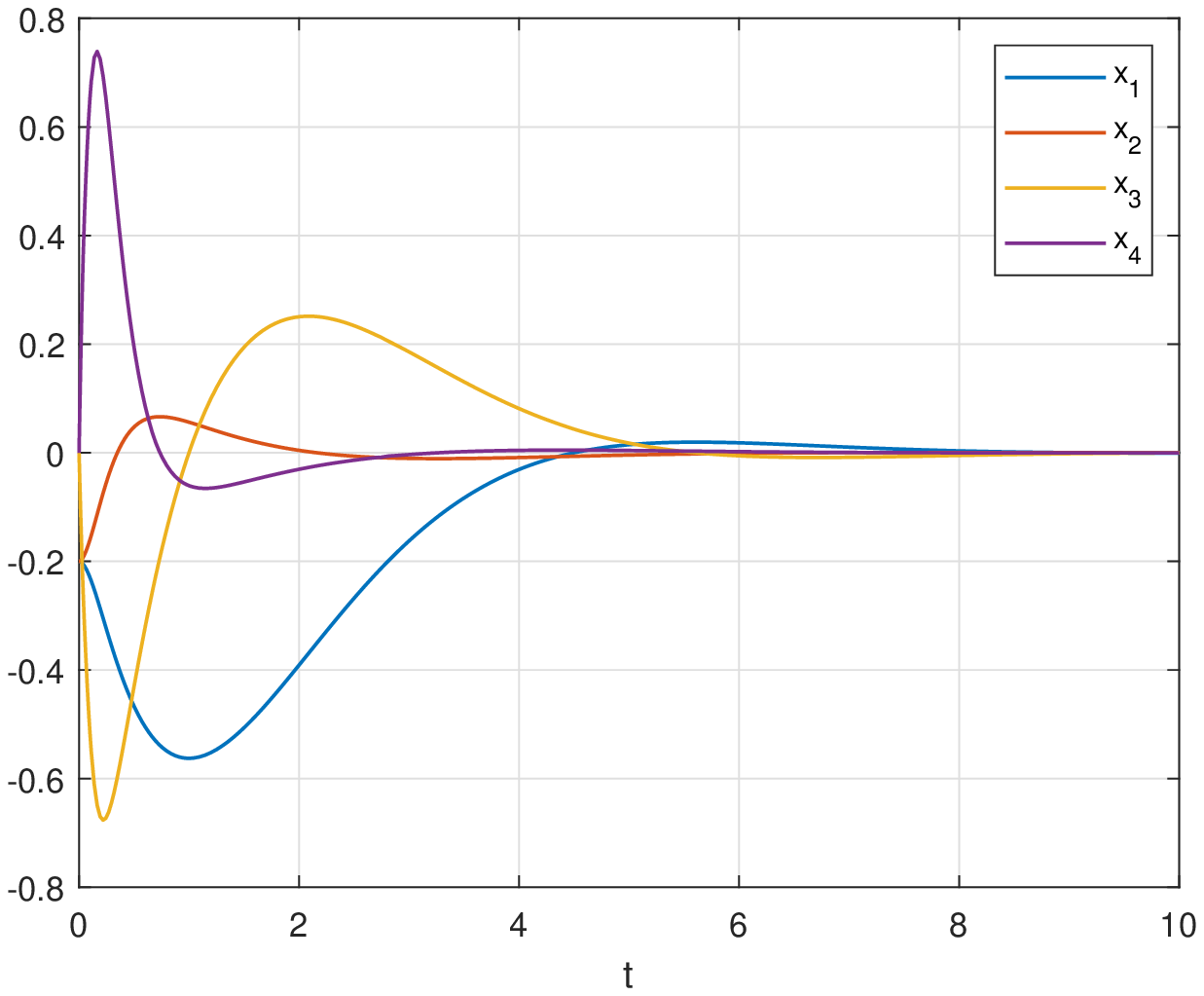} 
\includegraphics[scale=0.4]{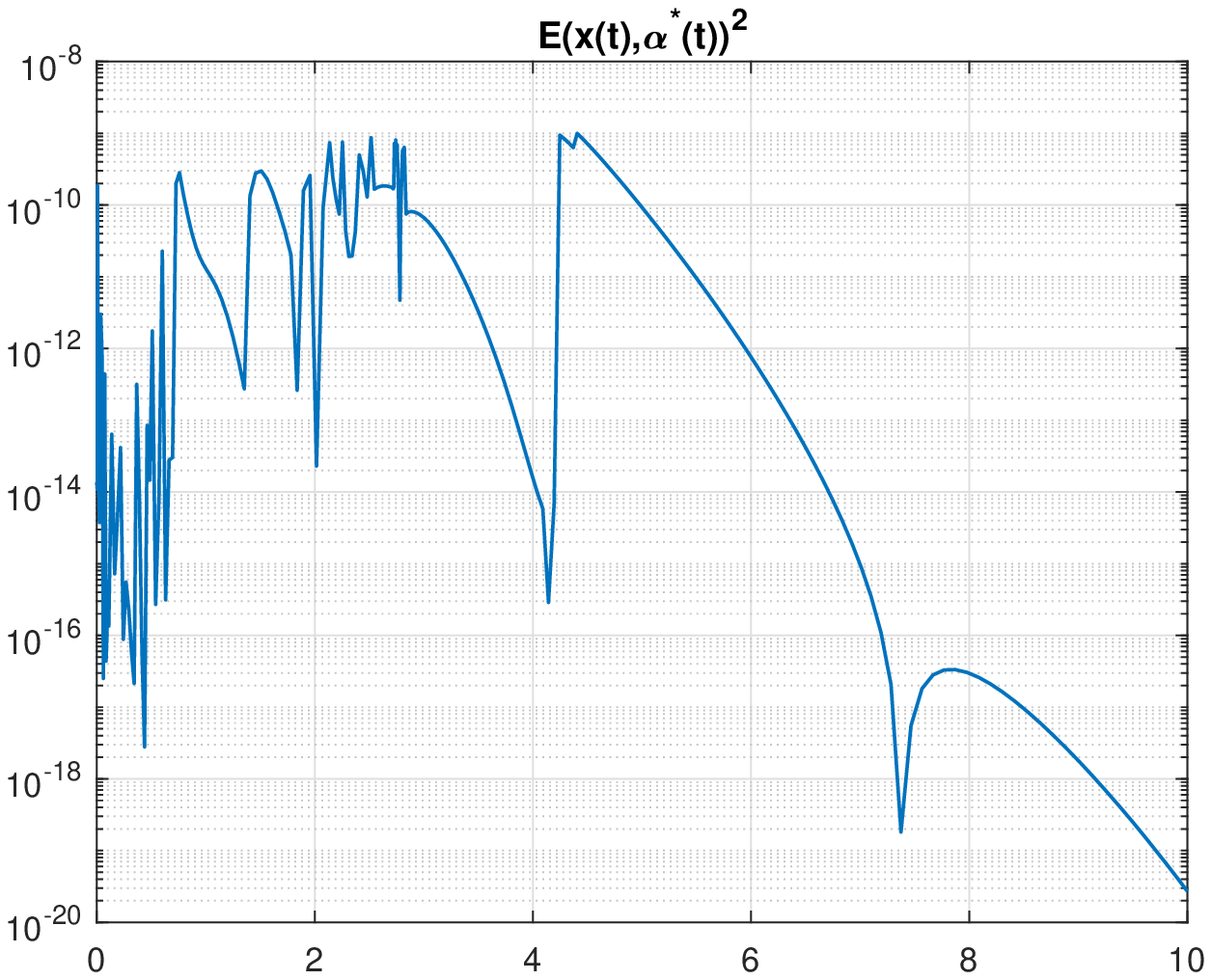}
                \caption{Test 2:  Optimal trajectories with initial condition $(-0.2,-0.2,0,0)$ (top) and the residual along the trajectory (bottom) with the optimal $\alpha^*$ .}
       \label{fig:case4}
	\end{figure}

Now let us fix the initial condition $(-0.2,-0.2,0,0)$. In this case the Riccati equations obtained using the semilinear representation formed by $A_0(x)$ can be solved on the optimal trajectory, allowing a comparison with the method based on the optimal affine combination. Again we notice that the introduction of the affine combination leads to a lower total residual and also to a lower total cost functional, as it can be observed by Table \ref{Table_err2}. The computation of the affine combination requires the resolution of minimization problems which leads to a slow-down of almost $8$ times. The optimal trajectory and the residual $E(x(t),\alpha(t))$ picking the optimal $\alpha^*$ are displayed in Figure \ref{fig:case4}. The residual presents a chattering behaviour at the beginning, due to the increase of the last two variables in the corresponding time interval. Finally we note again that as the dynamics approaches the equilibrium, the residual decays.

	\begin{table}[htbp]

		\centering
		\begin{tabular}{c|ccc}
					\toprule
					Semilinear form &Total cost & Total res & CPU   \\
					\midrule
$A_0(x)$ & 1.29 &  0.25 & \textbf{0.3s}\\
$\mathcal{A}(x,\alpha^*)$ & \textbf{1.27} &   \textbf{7.7e-10} & 2.5s\\

			\bottomrule		
		\end{tabular}	
	\caption{Test 2: Comparison between the fixed and the optimal case in terms of total cost and total residual. \label{Table_err2}}
	\end{table}	
	
\section{Conclusions}

We presented a novel method for the resolution of optimal control problems via a SDRE approach. The choice of a fixed semilinear representation for the entire resolution of the optimal control problem using the SDRE approach may be much less accurate than the application of the original HJB equations and sometimes it may affect the asymptotic stabilization of the dynamics. In this work we proposed a systematic construction of semilinear representations obtained perturbing the entries of an initial matrix and minimizing the residual in approximating the HJB by the SDRE. The numerical tests demonstrated that the proposed method is able to achieve low orders residual and a better accuracy in terms of cost functional. Moreover \ls{in the considered tests} it succeeds in finding a matrix $\mathcal{A}(x,\alpha^*)$ such that the couple $(\mathcal{A}(x,\alpha^*), g(x))$ turns to be stabilizable. \ls{The choice of the optimal semilinear form would be as expensive as solving the original HJB equation. Here we propose an approach which scales polynomially in the dimension, hence it is feasible from a numerical point of view.} Although the entire procedure has been presented as an online phase, an offline computation of the value function on a grid can be introduced. \ls{For instance, one may consider a data-driven approach for the construction of the value function (\cite{ABK21,dolgov2022data}).} In the future we aim at further investigating the combination of the semilinearizations and their minimization to achieve better results and their use for higher dimensional applications.

\bibliography{ifacconf}             % bib file to produce the bibliography

\begin{thebibliography}{20}
\providecommand{\natexlab}[1]{#1}
\providecommand{\url}[1]{\texttt{#1}}
\providecommand{\urlprefix}{URL }
\expandafter\ifx\csname urlstyle\endcsname\relax
  \providecommand{\doi}[1]{doi:\discretionary{}{}{}#1}\else
  \providecommand{\doi}{doi:\discretionary{}{}{}\begingroup
  \urlstyle{rm}\Url}\fi

\bibitem[{Albi et~al.(2021)Albi, Bicego, and Kalise}]{ABK21}
Albi, G., Bicego, S., and Kalise, D. (2021).
\newblock Gradient-augmented supervised learning of optimal feedback laws using
  state-dependent riccati equations.
\newblock \emph{IEEE Control Systems Letters}, 6, 836--841.

\bibitem[{Alla and Saluzzi(2020)}]{AS20}
Alla, A. and Saluzzi, L. (2020).
\newblock A {HJB}-{POD} approach for the control of nonlinear {PDEs} on a tree
  structure.
\newblock \emph{Applied Numerical Mathematics}, 155, 192--207.

\bibitem[{Azmi et~al.(2021)Azmi, Kalise, and Kunisch}]{AKK21}
Azmi, B., Kalise, D., and Kunisch, K. (2021).
\newblock Optimal feedback law recovery by gradient-augmented sparse polynomial
  regression.
\newblock \emph{J. Machin. Learn. Res.}, 22(48), 1--32.

\bibitem[{Banks et~al.(2007)Banks, Lewis, and Tran}]{BLT07}
Banks, H.T., Lewis, B.M., and Tran, H.T. (2007).
\newblock Nonlinear feedback controllers and compensators: a state-dependent
  riccati equation approach.
\newblock \emph{Comput. Optim. Appl.}, 37(2), 177--218.

\bibitem[{Beeler et~al.(2000)Beeler, Tran, and Banks}]{BTB00}
Beeler, S.C., Tran, H.T., and Banks, H.T. (2000).
\newblock Feedback control methodologies for nonlinear systems.
\newblock \emph{J. Optim. Theory Appl.}, 107(1), 1--33.

\bibitem[{Benner and Heiland(2018)}]{heiland}
Benner, P. and Heiland, J. (2018).
\newblock Exponential stability and stabilization of extended linearizations
  via continuous updates of riccati-based feedback.
\newblock \emph{International Journal of Robust and Nonlinear Control}, 28(4),
  1218--1232.

\bibitem[{Chow et~al.(2019)Chow, Darbon, Osher, and Yin}]{chow1}
Chow, Y.T., Darbon, J., Osher, S., and Yin, W. (2019).
\newblock Algorithm for overcoming the curse of dimensionality for
  state-dependent {H}amilton-{J}acobi equations.
\newblock \emph{J. Comput. Phys.}, 387, 376--409.

\bibitem[{Cloutier et~al.(1997{\natexlab{a}})Cloutier, D'Souza, and
  Mracek}]{sdrefirst1}
Cloutier, J.R., D'Souza, C.N., and Mracek, C.P. (1997{\natexlab{a}}).
\newblock Nonlinear regulation and nonlinear {$H_\infty$} control via the
  state-dependent {R}iccati equation technique. {I}. {T}heory.
\newblock In \emph{First {I}nternational {C}onference on {N}onlinear {P}roblems
  in {A}viation and {A}erospace ({D}aytona {B}each, {FL}, 1996)}, 117--130.
  Embry-Riddle Aeronaut. Univ. Press, Daytona Beach, FL.

\bibitem[{Cloutier et~al.(1997{\natexlab{b}})Cloutier, D'Souza, and
  Mracek}]{sdrefirst2}
Cloutier, J.R., D'Souza, C.N., and Mracek, C.P. (1997{\natexlab{b}}).
\newblock Nonlinear regulation and nonlinear {$H_\infty$} control via the
  state-dependent {R}iccati equation technique. {II}. {E}xamples.
\newblock In \emph{First {I}nternational {C}onference on {N}onlinear {P}roblems
  in {A}viation and {A}erospace ({D}aytona {B}each, {FL}, 1996)}, 131--141.
  Embry-Riddle Aeronaut. Univ. Press, Daytona Beach, FL.

\bibitem[{Dolgov et~al.(2021)Dolgov, Kalise, and Kunisch}]{tensor4}
Dolgov, S., Kalise, D., and Kunisch, K. (2021).
\newblock {T}ensor {D}ecompositions for {H}igh-dimensional
  {H}amilton-{J}acobi-{B}ellman {E}quations.
\newblock \emph{SIAM J. Sci. Comput.}, 43, A1625--A1650.

\bibitem[{Dolgov et~al.(2022)Dolgov, Kalise, and Saluzzi}]{dolgov2022data}
Dolgov, S., Kalise, D., and Saluzzi, L. (2022).
\newblock Data-driven {Tensor Train Gradient Cross Approximation for
  Hamilton-Jacobi-Bellman Equations}.
\newblock \emph{arXiv preprint arXiv:2205.05109}.

\bibitem[{Gr\"{u}ne and Pannek(2011)}]{larsnmpc}
Gr\"{u}ne, L. and Pannek, J. (2011).
\newblock \emph{Nonlinear model predictive control}.
\newblock Communications and Control Engineering Series. Springer, London.
\newblock Theory and algorithms.

\bibitem[{Gr\"{u}ne and Rantzer(2008)}]{GR2008}
Gr\"{u}ne, L. and Rantzer, A. (2008).
\newblock On the infinite horizon performance of receding horizon controllers.
\newblock \emph{IEEE Trans. Automat. Control}, 53(9), 2100--2111.

\bibitem[{Han et~al.(2018)Han, Jentzen, and E}]{ml1}
Han, J., Jentzen, A., and E, W. (2018).
\newblock Solving high-dimensional partial differential equations using deep
  learning.
\newblock \emph{Proc. Natl. Acad. Sci. USA}, 115(34), 8505--8510.

\bibitem[{Herty and Kalise(2018)}]{HK18}
Herty, M. and Kalise, D. (2018).
\newblock Suboptimal nonlinear feedback control laws for collective dynamics.
\newblock In \emph{2018 IEEE 14th International Conference on Control and
  Automation (ICCA)}, 556--561.

\bibitem[{{Jones} and {Astolfi}(2020)}]{Astolfi2020}
{Jones}, A. and {Astolfi}, A. (2020).
\newblock On the solution of optimal control problems using parameterized
  state-dependent {R}iccati equations.
\newblock In \emph{2020 59th IEEE Conference on Decision and Control (CDC)},
  1098--1103.

\bibitem[{Krener(2020)}]{krener2020brekht}
Krener, A.J. (2020).
\newblock Al’brekht’s method in infinite dimensions.
\newblock In \emph{2020 59th IEEE Conference on Decision and Control (CDC)},
  5653--5658. IEEE.

\bibitem[{Nakamura-Zimmerer et~al.(2021)Nakamura-Zimmerer, Gong, and
  Kang}]{nakazim}
Nakamura-Zimmerer, T., Gong, Q., and Kang, W. (2021).
\newblock {A}daptive {D}eep {L}earning for {H}igh-{D}imensional
  {H}amilton--{J}acobi--{B}ellman {E}quations.
\newblock \emph{SIAM J. Sci. Comput.}, 43(2), A1221--A1247.

\bibitem[{Wang and Wu(1998)}]{WANG98}
Wang, J. and Wu, G. (1998).
\newblock A multilayer recurrent neural network for solving continuous-time
  algebraic {R}iccati equations.
\newblock \emph{Neural Networks}, 11(5), 939--950.

\bibitem[{Çimen(2008)}]{CIMEN}
Çimen, T. (2008).
\newblock State-dependent riccati equation (sdre) control: A survey.
\newblock \emph{IFAC Proceedings Volumes}, 41(2), 3761--3775.
\newblock \doi{https://doi.org/10.3182/20080706-5-KR-1001.00635}.
\newblock 17th IFAC World Congress.

\end{thebibliography}
                                                     % with bibtex (preferred)
                                                   
%\begin{thebibliography}{xx}  % you can also add the bibliography by hand

%\bibitem[Able(1956)]{Abl:56}
%B.C. Able.
%\newblock Nucleic acid content of microscope.
%\newblock \emph{Nature}, 135:\penalty0 7--9, 1956.

%\bibitem[Able et~al.(1954)Able, Tagg, and Rush]{AbTaRu:54}
%B.C. Able, R.A. Tagg, and M.~Rush.
%\newblock Enzyme-catalyzed cellular transanimations.
%\newblock In A.F. Round, editor, \emph{Advances in Enzymology}, volume~2, pages
%  125--247. Academic Press, New York, 3rd edition, 1954.

%\bibitem[Keohane(1958)]{Keo:58}
%\newblock \emph{Power and Interdependence: World Politics in Transitions}.
%\newblock Little, Brown \& Co., Boston, 1958.

%\bibitem[Powers(1985)]{Pow:85}
%\newblock Is there a way out?
%\newblock \emph{Harpers}, pages 35--47, June 1985.

%\bibitem[Soukhanov(1992)]{Heritage:92}
%\newblock \emph{{The American Heritage. Dictionary of the American Language}}.
%\newblock Houghton Mifflin Company, 1992.

%\end{thebibliography}
           % Sections and subsections are supported  
                                                                         % in the appendices.
\end{document}